\documentclass[a4paper,preprint]{imsart}
\usepackage[authoryear]{natbib}
\usepackage[english]{babel}
\usepackage{amsmath,amsfonts,amsthm, amssymb}
\usepackage{graphicx}
\usepackage{enumerate}
\usepackage{color}
\usepackage{subfigure}
\usepackage{hyperref}
\usepackage{framed}

\newcommand{\E}{\mathbb{E}}

\newcommand{\1}{\textbf{1}}

\newcommand{\R}{\mathbb{R}}
\renewcommand{\P}{\mathbb{P}}
\newcommand{\Xn}{\mathbf{X}^n}
\newcommand{\hp}{\widehat{p}}

\newtheorem{cond}{Condition}

 \newtheorem{assumption}{Assumption}

 \theoremstyle{plain}
 \newtheorem{thm}{Theorem}[section]
 \newtheorem{lem}[thm]{Lemma}

 \theoremstyle{definition}
 \newtheorem{defn}[thm]{Definition}

 \theoremstyle{remark}

\begin{document}

\begin{frontmatter}

\title{Risk quantification for the thresholding rule for multiple testing using Gaussian scale mixtures}
\runtitle{Testing risks for Gaussian scale mixtures priors}
\begin{aug}
\author{\fnms{Jean-Bernard}  \snm{Salomond}\ead[label=e1]{jean-bernard.salomond@u-pec.fr}}
  \address{Universit\'e Paris-Est, Laboratoire d'Analyse et de Math\'ematiques Appliqu\'ees (UMR 8050)
UPEM, UPEC, CNRS, F-94010, Cr\'eteil, France\\ 
          \printead{e1}}
 \runauthor{J-B Salomond}
\end{aug}

\begin{abstract}{
In this paper we study the asymptotic properties of Bayesian multiple testing procedures for a large class of Gaussian scale mixture priors. We study two types of multiple testing risks: a Bayesian risk proposed in \cite{MR2850212} where the data are assume to come from a mixture of normal, and a frequentist risk similar to the one proposed by \cite{arias-castro2017}. Following the work of \cite{vanderpas2016}, we give general conditions on the prior such that both risks can be bounded. For the Bayesian risk, the bound is almost sharp. This result show that under these conditions, the considered class of continuous prior can be competitive with the usual two-group model (e.g. spike and slab priors). We also show that if the non-zeros component of the parameter are large enough, the minimax risk can be made asymptotically null. The separation rates obtained are consistent with the one that could be guessed from the existing literature \citep[see][]{vanderpas2017UQHS}. For both problems, we then give conditions under which an adaptive version of the result can be obtained. 
}
\end{abstract}



\end{frontmatter}

\maketitle

\section{Introduction}

Multiple testing has become a topic of particular interest over the past decades. High dimensional models are now quite common for instance in genomics, bio-informatic or even finance. In all theses applications, it is now well known that multiple testing has to be accounted for \citep[see for instance][]{MR1946571,efron2002empirical,MR2065197,MR2109489,MR2235051}. 
For high dimensional models, \cite{MR2980074} and \cite{rossell2015non} proposed methods based on non local prior.
Recently, the multiple testing problem has been particularly well studied when the number of true positive is supposed to be small with respect to the total amount of tests. This problem is closely related to estimation under sparsity assumption which is now a well known problem. Frequentists approaches such as the LASSO \cite{Tibshirani1996} are now well understood, and several Bayesian approaches have been developed. In particular \cite{Castillo2012,Castillo2015} studied in detail the spike and slab prior and obtained the minimax posterior contraction rate for the posterior, and some desirable properties for recovering the true model.
Another celebrated class of prior are the so called \emph{one class} priors (or one group priors). Among others, the Horseshoe prior has been studied in \cite{Carvalho2010} and \cite{vanderPas2014}, \cite{Ghosh2015} studied a generalization of the Horseshoe type priors, \cite{Rockova2015} proposed a continuous version of the spike and slab prior. Recently \cite{vanderpas2016} proposed some general conditions on the one group prior to obtain minimax posterior concentration rate. However these results are focused on estimation of the parameter rather than testing.



We consider the idealistic Gaussian sequence model
\begin{equation}
X_i = \theta_i + \epsilon_i , ~ i= 1, \dots,n ,%
\label{eq:model}%
\end{equation}
where the $\epsilon_i$ are independent and identically distributed $\mathcal{N}(0,v^2)$. In the sequel we should fix $v^2$ to be $1$ for simplicity.
We then consider a Gaussian scale mixture prior, where each coefficient $\theta_i$ receives a prior
\begin{equation}
\theta_i | \sigma_i^2 \sim \mathcal{N}(0,\sigma_i^2), ~ \sigma_i^2 \sim \pi(\sigma^2_i),
\label{eq:prior}
\end{equation}

where $\pi$ is a density on the positive reals. Here the random parameter $\sigma_i^2$ both accounts for the overall sparsity of the parameter, but is still flexible enough to detect large signals. 
This class of prior is general enough to encompass most of the one group priors introduced in the literature so far, and all those stated above in particular.
The theoretical properties of these priors are now well known for estimation, and they are widely used in practice.
For these priors we can easily derive the posterior distribution on the $\theta_i$ and $\sigma_i^2$
\begin{align*}
\pi(\sigma_i|X_i) \propto (1+\sigma_i)^{-1/2} e^{-\frac{X_i^2}{2} \frac{\sigma_i^2}{1+\sigma_i^2}}\pi(\sigma_i), \\ \theta_i|\sigma_i^2,X_i \sim \mathcal{N}\left(X_i\frac{\sigma_i^2}{1+\sigma_i^2}, \frac{\sigma_i^2}{1+\sigma_i^2}\right).
\end{align*}
The parameter $\kappa_i = \frac{\sigma_i^2}{1+\sigma_i^2}$ is often call a shrinkage parameter and play a crucial role in both the estimation and testing for these models.

In this paper we are interested in the problem of testing simultaneously multiple hypotheses on the parameter $\theta$. More precisely, we consider the following sequence of test 
$$
H_{0,i} : \theta_i = 0 \text{ versus } H_{1,i}: \theta \neq 0,~  i=1,\dots, n.
$$
Denote by $S_0$ the unknown support of true \emph{signals} or true rejections i.e. $S_0 = \{ i, \theta_i \neq 0 \}$, and $S_0^c$ the support of the true positive $S_0^c = \{i, \theta_i = 0\}$. 
Consider a sequence of tests $\Xi = (\xi_1, \dots \xi_n)$. A test $\xi_i$ for hypothesis $H_{0,i}$ versus $H_{1,i}$ is a statistic of the observed data $\Xn=(X_1, \dots, X_n)$ such that $\xi_i = 0$ if $H_{0,i}$ is accepted and $\xi_i = 1$ otherwise. 
A multiple testing procedure returns a subset of $\{1, \dots, n\}$ of hypotheses that are accepted, representing the support of the \emph{signals}. This is equivalent to returning the support of the vector $\Xi.$

A great advantage of the spike and slab approach for sparse models is that it provides exact $0$ under the posterior, and one can thus easily derive probability of inclusion and do multiple testing. 
One group prior models that are considered here, since they are continuous, do not share this property. Nonetheless in \cite{Polson2010}, by analogy with the Spike and Slab approach propose to call \emph{signals} parameters $\theta_i$ for which $\E(\kappa_i|X_i)\geq 1/2$. The authors have investigated this thresholding rule for the Horseshoe prior and observed empirically that it behaved similarly to the Spike and Slab inclusion probability. Theoretical properties of this approach to multiple testing were later studied in \cite{Ghosh2015} for some specific priors and a additive loss function derived from \cite{MR2850212}. In the latter paper the authors derived the optimal risk for this additive loss. 

Here we will study optimality of such multiple testing rule when only a small number $p_n \ll n$ of the parameters are true signals for different risks.
In the same spirit of \cite{vanderpas2016} that gave conditions on the prior to bound the posterior contraction rate and convergence rate of a Bayesian estimator, in this paper, we give conditions on the prior distribution under which the we can control the Type I and Type II error for each individual test. This allows us to give upper bounds on the Bayes Risk as studied in \cite{Ghosh2015}. 
The additive Bayes Risk studied in \cite{MR2850212} and \cite{Ghosh2015} is a Bayesian integrated risk for a two group prior 
$$
\mu: \theta \sim (1-p)\delta_{0} + p \mathcal{N}(0,\psi^2),
$$
which in turns give the following marginal distribution for the data 
$$
X_i \sim (1-p)\mathcal{N}(0,1) + p\mathcal{N}(0,\psi^2).
$$
For a sequence of tests $\Xi = (\xi_1, \dots, \xi_n)$, \cite{MR2850212}, following ideas that goes back to \cite{MR0084952}, propose an overall additive risk
$$
R^\mu(\Xi) = \sum_{i=1}^{n} (1-p)\P_{\mathcal{N}(0,1)}(\xi_i = 1) + p\P_{\mathcal{N}(0,\psi^2)}(\xi_i = 0). 
$$
Similarly to the proposition of \cite{Carvalho2010} to take $\E(\kappa_i|X_i)$ as a proxy of the inclusion probability, we will consider tests of the form $\xi_i = \1\{m_{X_i} > \alpha\}$ where $m_x = \E(\kappa_i|X=x)$ and $\alpha$ is some fixed threshold. Thus taking $\alpha=1/2$ will result in the same procedure as the one studied in \cite{Carvalho2010} and \cite{Ghosh2015}.

We also consider a minimax approach to the multiple testing problem. Minimax theory for testing has been introduced by \cite{Ingster93}, however, the setting considered in this paper does not seems suited here. We focus on another optimality criterion based on separation rate of the test, that is the rate at which parameters in the alternative can approach the null hypothesis and still be detected by the test. In single hypothesis testing, separation rates have been studied in \cite{baraud2002}. In a Bayesian setting \cite{salomondTest} proposed general way of constructing tests and control their separation rate. Recently \cite{fromont2016} proposed a version of the separation rate for testing in a multiple testing setting. In a similar spirit \cite{arias-castro2017} studied a multiple testing risk based on the False Discovery Rate and False Nondiscovery Rates. 
False Discovery Rate (FDR) is a measure of the Type I error in a multiple testing setting. Let $\Xi = (\xi_1, \dots, \xi_n)$ be a sequence of tests and let $S_0$ be the sets of \emph{true} alternative hypotheses, the False Discovery Rate is defined as 
$$
FDR_n(\Xi) = \E\left( \frac{\sum_{i\in S_0^c} \xi_i }{\sum_{i=1}^n \xi_i \vee 1 } \right).
$$
Since the seminal work \cite{benjamini1995controlling} there has been an abundant literature on how to control the FDR for different types of models. 
However only a few results are concerned with the control of the power of such procedure. Following \cite{arias-castro2017} we study the False Nondiscovery Rate (FNR) which is the expectation of the proportion of false negative test, defined as 
$$
FNR_n(\Xi) = \E \left( \frac{\sum_{i\in S_0}(1- \xi_i) }{p_n} \right). 
$$

We consider the same risk as in \cite{arias-castro2017} or \cite{Rabinovich17}
\begin{equation}
R^{sup}_n(\Xi) = \{FDR_n(\Xi) + FNR_n(\Xi)\}. 
\label{eq:intro:risk:sup}
\end{equation}
In particular we will aim at bounding the supremum of this risk over a set of parameters
$$
\sup_{\theta, \theta \in T_n^c } R_n^{sup} (\Xi), 
$$
where $T_n$ is some neighborhood of $0$. Intuitively $T_n$ is the set of parameters that cannot be detected by the test. We chose $T_n$ of the form $T_n:= \{ \theta, d_n(\theta,0)\leq \rho_n \}$ and thus $\rho_n$ is the detection boundary. The smallest the $\rho_n$ is, the sharper the test. In this paper, we will show that under some conditions on the prior, the considered test will have asymptotically null risk if $\rho_n \gtrsim \sqrt{2\log(n/p_n)}$ for $d_n$ the $l_\infty$ distance. It is well known that for our model, parameters bellow $\sqrt{2\log(n)}$ will be shrunk toward $0$, implying a consequent bias for estimation. However it seems that parameter between $\sqrt{2\log(n/p_n)}$ and $\sqrt{2\log(n)}$ can still be detected by the test.  

We also study empirical Bayes approaches to the testing problem to adapt to the more realistic case where the number of signals $p_n$ is unknown, and give some conditions on the estimator $\hat{p}$ of $p_n$ such that both risks can be bounded. Other types of adaptive procedures have been proposed in the literature, such as empirical Bayes based on the Marginal Maximum likelihood estimator or fully Bayes procedure for the Horseshoe prior \citep[see][]{vanderpas2017}. However, these procedures are difficult to adapt to the more general set up we consider here, and will thus not be treated.  

The paper is organized as follows, in section \ref{sec:main:result} we state the main results providing conditions on the prior such that we can bound the type I and type II error for each individual test. We then give upper bounds for the Bayes Risk and give an upper bound on the posterior separation rate for the multiple testing risk. We then give the adaptive version of the main results in section \ref{sec:adaptive}. Section \ref{sec:proofs} is devoted to the proofs. 

\paragraph{Notations} Throughout the rest of the paper we will denote by $a\wedge b$ and $a \vee b$ respectively the minimum and the maximum between to real numbers $a$ and $b$. The probability density function of a standard normal will be denoted $\phi$ and its cumulative distribution function $\Phi$. We will denote $a\lesssim b$ if there exist an absolute constant $C$ such that $a \leq C b$ and $a \asymp b$ if $a\lesssim b $ and $b\lesssim a$. We will denote by $\tau_n(p) = p/n$ and $\nu_n(p) = \sqrt{\log(n/p)}$.

\section{Prior conditions}
\label{sec:main:result}

Because scale mixture priors are continuous, we do not have access to inclusion probability to perform multiple testing.
Because the $\theta_i$ are a posteriori independent, we can use the coordinatewise posterior. In particular, we get that 
\begin{align}
\E(\kappa_i|X_i = x) &= { \int_0^1 z (1-z)^{-3/2} e^{\frac{x^2}{2}z} \pi\left(\frac{z}{1-z}\right) dz \over \int_0^1 (1-z)^{-3/2} e^{\frac{x^2}{2}z} \pi\left(\frac{z}{1-z}\right) dz } \nonumber  \\
&= { \int_0^\infty u(1+u)^{-3/2} e^{\frac{x^2}{2} \frac{u}{1+u}} \pi(u) du \over \int_0^\infty (1+u)^{-1/2} e^{\frac{x^2}{2} \frac{u}{1+u}} \pi(u) du}. \label{eq:Ekappa}
\end{align}
From Tweedie's formula \citep{Robbins1956}, we have that the posterior mean of $\theta_i$ given an observation $x_i$ is equal to $x_i + \frac{d}{dx}p(x_i)$ where $p(x_i)$ is the marginal distribution of $x_i$. This in turns gives that the posterior mean of $\theta_i$ given $X_i$ is $\widehat{\theta}_i = \E(\kappa_i|X_i) X_i$. An advantage of scale mixtures over spike-and-slab types of priors is that the integrals in \eqref{eq:Ekappa} can be computed via integral approximation methods see \citep[see][in the context of the horseshoe prior]{vanderPas2014,Carvalho2010}. This is also a great advantage for multiple testing as $\E(\kappa_i|X_i)$ plays a central role in the approach proposed in \cite{Carvalho2010}. In particular, because the parameters are a posteriori independent, we will be able to bound each individual type I and type II error. 

Our main results are conditions on the prior under which we can guaranty upper bounds on multiple testing risks. The conditions are similar to the ones proposed in \cite{vanderpas2016} under which the authors prove some upper bounds on the contraction rates of the order of the minimax rate. We first state the conditions. Condition \ref{cond1} is required to bound the type II error of each individual tests. The remaining two are used to bound each individual type I error. 

The first condition involves a class of uniformly regularly varying function at infinity. defined as follows. 
\begin{defn}
A function $L$ is said to be uniformly regularly varying at infinity if there exists $R, u_0 > 0$ such that 
$$
\forall a \in [1,2],~ \forall u\geq u_0, ~  \frac{1}{R} \leq \frac{L(au)}{L(u)} \leq R.
$$ 
\end{defn} 
In particular, polynomial functions $L(u) = u^b$ and powers of logarithm $L(u) = \log^b(u)$ with $b\in \R$ are uniformly regular varying at infinity. For more details on uniformly regular varying function, see \cite{vanderpas2016}. The only difference with their definition is that here we allow $u_0$ to be less than 1. We now present the first condition. 

\begin{cond}\label{cond1}
For some $b\geq 0,$ we can write $u\mapsto \pi(u) = L_n(u) e^{-bu},$ where $L_n$ is a uniformly regularly varying function at infinity for some $R, u_0> 0$ which do not depend on $n.$ Suppose further that there are constants {$C', b'>0,$} $K\geq 0,$ and $u_* \geq 1,$ such that 
\begin{align}
	  C' \pi(u) \geq \tau_n(p)^{K} e^{-b'u}  \quad \text{for all} \ u\geq u_*.
	\label{eq.assump_on_lb_Ln}
\end{align}
\end{cond}

Condition \ref{cond1} is a sufficient condition on the tails of the prior to bound the type II error of each individual test. This condition states that the tails of $\pi$ should decay at most exponentially fast, up to a uniformly regularly varying function. This is consistent with the conditions on the slab distribution in \cite{Castillo2012}. The following two conditions unsure that the prior penalizes enough small signals, which will be necessary to bound the individual type I errors. Condition \ref{cond2} requires that $\pi$ puts enough mass on a neighborhood of $0$ while Condition \ref{cond3} describes the decay of $\pi$ away from $0$. 

\begin{cond}\label{cond2}
Suppose that there is a constant $c>0$ such that $\int_0^1 \pi(u) du \geq c.$
\end{cond}
Define $s_n$ as 
\begin{align}
	s_n:= \tau_n(p) \nu_n(p)^2,
	\label{eq.sn_def}
\end{align}
then to get bound on the different risks, we ask that $\pi$ satisfies the following condition. 
\begin{cond}[$C$] \label{cond3}
Assume that there is a constant $C,$ such that
\begin{align*}
	\int_{s_n}^\infty \Big( u \wedge \frac{\nu_n(p)^3}{\sqrt{u}} \Big) \pi(u) du
	+ \nu_n(p) \int_1^{\nu_n(p)^2} \frac{\pi(u)}{\sqrt{u}} du
	\leq C s_n.
\end{align*}
\end{cond}

We let Condition \ref{cond3} depend on the constant $C$ has it will appears in the bound on the different risks. In particular, it will appears that by tuning $C$ appropriately or by taking $C = C_n \to 0$ with $n$, we can obtain almost sharp bounds especially for the Bayes risk. 
%

\section{Deterministic $p$}
In this section we will assume that $p_n$ is known. 
Following the results of \cite{MR2850212}, we assume that the prior $\mu$ considered for the Bayes risk $R^\mu$ satisfies the following assumption. 

\begin{assumption} The distribution $\mu$ is of the form $$\mu = (1-\frac{p_n}{n}) \mathcal{N}(0,1) + \frac{p_n}{n}\mathcal{N}(0,\psi_n),$$ 
\label{ass:mu}
 with ${p_n}/{n} = o(1)$ and $\psi_n = C_\psi^{-1} \log(n/p_n)(1+o(1))$.
\end{assumption}
Under this assumption, \cite{MR2850212} showed that the Bayes Oracle $\Xi^*$ as the asymptotic risk 
\begin{equation}
R^\mu(\Xi^*) = p_n(2\Phi(\sqrt{C_\psi}) - 1)(1 + o(1)).
\label{eq:RiskOracle}
\end{equation}
In their work \cite{Ghosh2015} proved that for certain class of priors, the optimal risk \eqref{eq:RiskOracle} could be obtained. The following theorem gives an upper bound on the risk that depends on the constants in Conditions \ref{cond1}-\ref{cond3}. 

\begin{thm}
\label{thm:BayesRisk}
Work under model \eqref{eq:model} and assume that the prior is of the form \eqref{eq:prior} and satisfies Conditions \ref{cond1}-\ref{cond3}($C$) with $p = p_n$. Assume that the distribution $\mu$ satisfies the Assumption \ref{ass:mu}. Consider for the sequence of test $\xi_i = \1\left\{ \E(\kappa_i|X_i) > \alpha \right\}$ for any fixed $\alpha \in (0,1)$, then
\begin{equation}
R^\mu(\xi_1, \dots \xi_n) \leq p_n\left( \frac{8\sqrt{\pi} C}{c \alpha} + 2\Phi\left( \sqrt{2K(u_0+1)C_\psi}\right)-1 \right)(1+o(1)).
\label{eq:bound:thm:BayesRisk}
\end{equation}
\end{thm}
The proof of Theorem \ref{thm:BayesRisk} is postpone to section \ref{sec:proofs}. 
This theorem thus gives bounds on the Bayes Risk for a wide class of priors, but is slightly less sharp than the ones obtained in \cite{Ghosh2015} for instance. In particular, when considering the horseshoe type priors as in \cite{Ghosh2015} with $\alpha_n= 1/2$, the bound \eqref{eq:bound:thm:BayesRisk} becomes
$$
R^\mu(\Xi) \leq p_n\left( 2\Phi\left( \sqrt{2C_\psi}\right)-1 \right)(1+o(1)),
$$
to be compared with the oracle bound from \eqref{eq:RiskOracle}. The suboptimal constant is an artifact of the proof of the general case. More precisely, because we do not impose any particular form for the prior $\pi$ on the near $0$, one cannot control precisely the Type II error. 

We now go on and study the minimax Risk. Under the same conditions, the following theorem gives an upper bound on the separation such that the minimax risk \eqref{eq:intro:risk:sup} can be bounded. 
\begin{thm}
\label{thm:seprate}
Work under model \eqref{eq:model} and assume that the prior is of the form \eqref{eq:prior} and satisfies Conditions \ref{cond1}-\ref{cond3}($C$) with $p = p_n$. For all $v_n\to \infty$, let $\rho_n =C_1 +  \sqrt{2K(u_0+1)\log(n/p_n)} + v_n$ then we have, for all $\lambda \in (0,1)$
$$
\sup_{\theta, \min_{i\in S_0}|\theta_i| \geq \rho_n}R^{sup}_n(\Xi) \leq \left(\frac{1}{1+\frac{\lambda \alpha c}{8 C\sqrt{\pi}}} + \Phi(-v_n)\right)(1 + o(1))
$$
\end{thm}

The proof of theorem \ref{thm:seprate} is postpone to section \ref{sec:proofs}. We see that with a careful choice of the prior hyper-parameters, $\R^{sup}_n$ can be made asymptotically small.  
The separation rate of the test $\rho_n$ is of the order of $\sqrt{\log(n/p_n)}$. This is surprising as in estimation problems, the \emph{universal threshold} under which the parameter are shrunk toward $0$ is known to be $\sqrt{2\log(n)}$. However here it seems that this shrinkage is not too important for parameters between $\sqrt{2\log(n)}$ and $\rho_n$, which is essentially of the order of $\sqrt{2K(u_0+1)\log(n/p_n)}$, so that the proportion of false negative is not too high.

\section{Adaptive results}
\label{sec:adaptive}

In this section we give an adaptive counterpart to the results presented in section \ref{sec:main:result}. The conditions \ref{cond1}-\ref{cond3} depend on the number $p_n$ of \emph{true signals} in the parameter $\theta$. We give general conditions on an estimator $\hp$ of $p$ under which we can bound each risk. The assumptions are similar to the one used in \cite{vanderpas2017} for instance. We then show that this conditions are satisfied for a standard estimator of $p$. 

\subsection{Assumptions and main results}
We first give an assumption on $\hp$ in order to bound the Bayes risk $R^\mu_n$. Similarly to \cite{vanderpas2017,vanderpas2017UQHS}, we require that the estimator do not overestimate nor underestimate the number of true positive $p_n$. 
\begin{cond}
\label{ass:adapt:Bayes}
Let $\mu$ be some probability mesure satisfying assumption \ref{ass:mu}. For some absolute constants $C^u>0$, $c_d>0$ and $C_d\geq 0$ the estimator $\hat{p}$ of $p_n$ satisfies 
\begin{itemize}
\item There exists $C^u$ such that $\P^n_\mu(\hat{p} \leq C^u  p_n) =1 +  o(\frac{p_n}{n})$
\item There exists absolute constants $C_d\geq 0$, $c_d >0$ and $\zeta \geq 0$ such that
 $\P^n_\mu\left(\hp \geq c_d p_n \left(\frac{n}{p_n}\right)^{-\zeta}e^{-C_d\sqrt{K\log(n/p_n)}} \right) = 1 + o(1)$ 
\end{itemize} 
\end{cond}
The conditions is stronger here compare to the ones proposed in \cite{vanderpas2017,vanderpas2017UQHS}, as we require that the probability of overestimating $p_n$ is $o(p_n/n)$. This is due to the fact that when over estimating the number of true signals, the risk is of the order of $n$. Note that although this conditions is stronger than the one proposed in the literature, however the simple estimator of $p_n$ proposed in \cite{vanderPas2014} for instance will satisfy conditions \ref{ass:adapt:Bayes} under some additional assumptions on the number of true signals similar to the ones considered in \cite{Ghosh2015}. 
The second part of the condition requires that the estimator $\hp$ does not under estimate the number of true signal. Note that this condition is quite flexible and could me made trivial with a proper choice of $\zeta$ when $p_n$ is a of the polynomial order of $n$. However a larger $\zeta$ will deteriorate the bound on the risk. The following theorem gives a upper bound on the Bayes risk for the empirical Bayes method.

\begin{thm}
\label{thm:AdaptBayesRisk}
Work under model \eqref{eq:model} and assume that the prior is of the form \eqref{eq:prior} and satisfies Conditions \ref{cond1}-\ref{cond3}($C$) with $p = \hat{p}$ and let $\hat{p}$ be an estimator of $p_n$ satisfying Assumtion \ref{ass:adapt:Bayes}. Assume that the distribution $\mu$ satisfies the Condition \ref{ass:mu}. Consider for the sequence of test $\xi_i = \1\left\{ \E(\kappa_i|X_i) > \alpha \right\}$ for any fixed $\alpha \in (0,1)$, then
\begin{equation}
R^\mu(\xi_1^{\hat{p}}, \dots \xi_n^{\hat{p}}) \leq p_n\left( \frac{8\sqrt{\pi} C C^u}{c \alpha} + 2\Phi\left( \sqrt{2K(u_0+1)(1 + \zeta)C_\psi}\right)-1 \right)(1+o(1)).
\label{eq:bound:thm:AdaptBayesRisk}
\end{equation}
\end{thm}
The proof of theorem \ref{thm:AdaptBayesRisk} is postpone to appendix \ref{seq:proof:adaptive}. We see that if the Condition \ref{ass:adapt:Bayes} is satisfied with $\zeta=0$, we get a bound similar to the one obtained in the non-adaptive case. When $p_n$ is of a polynomial order of $n$, \cite{ghosh2016} have shown that the simple estimator considered in \cite{vanderPas2014}, satisfies condition \ref{ass:adapt:Bayes} with $\zeta =0$ (see their remark 4) . We also propose an adaptive version of Theorem \ref{thm:seprate} with however less strong conditions on the estimator of $p_n$.

\begin{cond}
\label{ass:adapt:Minimax}
For some sequence $\rho_n$, let $\theta$ be such that $\forall i \in S_0, |\theta_i| \geq \rho_n^0$. For some absolute constants $C^u>0$, $c_d>0$ and $C_d\geq 0$ the estimator $\hat{p}$ of $p_n$ satisfies 
\begin{equation}
\P_\theta^n \left( \gamma_n \leq \hp \leq  C^u  p_n)\right) =1 +  o(1)
\label{eq:condition:adapt:minimax}
\end{equation}
\end{cond}
This condition is fairly similar to Condition \ref{ass:adapt:Bayes}, however we only require that the probability of over estimating the number of signals goes to $0$. Here again by choosing $\nu$ large enough, the assumption on the lower bound on $\hp$ in \eqref{eq:condition:adapt:minimax} will be satisfied, however this will impact the separation rate of the test. 

\begin{thm}
\label{thm:adapt:seprate}
Work under model \eqref{eq:model} and assume that the prior is of the form \eqref{eq:prior} and satisfies Conditions \ref{cond1}-\ref{cond3}($C$) with $p = \hat{p}$ and let $\hat{p}$ be an estimator of $p_n$ satisfying Condition \ref{ass:adapt:Minimax} with $\rho_n^0$. For all $v_n\to \infty$, let $\rho_n =  C_1 +   \sqrt{2K(u_0+1)\log(n/\gamma_n)} + v_n$ then we have, for all $0 < \lambda < \Phi(v_n) $
$$
\sup_{\theta, \min_{i\in S_0}|\theta_i| \geq \rho_n \vee \rho_n^0} R^{sup}_n(\Xi) \leq \left(\frac{1}{1+\frac{\lambda \alpha c}{8 C^uC\sqrt{\pi}}} + \Phi(-v_n)\right)(1 + o(1))
$$
\end{thm}
The proof of this theorem is postponed to appendix \ref{seq:proof:adaptive}. The bounds are comparable with the non adaptive one. The separation rate is of the order of $\sqrt{\log(n/\gamma_n)}$ where $\gamma_n$ is a lower bound on the estimated number of true signals. Thus for the simple estimator, given the results of \cite{vanderpas2017,vanderpas2017UQHS} Condition \ref{ass:adapt:Minimax} holds true for $\gamma_n = 1$ and the separation rate is of the same order as the contraction rate. Moreover this rate is consistent with the ones that could be obtained using tests base on the confidence intervals obtained by \cite{vanderpas2017UQHS}.

\section{Proofs}
\label{sec:proofs}
\subsection{Proof of Theorem \ref{thm:BayesRisk}}
\subsubsection{Control of the type I Error} 

We first control the Type one error $\P_0(\E(\kappa_i|X_i) \geq \alpha )$, where $P_0$ is the probability distribution of a $\mathcal{N}(0,1)$. Following the proof of \cite{vanderpas2016} we have denote by $m_x = \E(\kappa_i|X_i = x)$, and we have 

\begin{equation}
m_x  = \frac{\int_0^1 z(1-z)^{-3/2} e^{\frac{x^2}{2}z} \pi\left( \frac{z}{1-z}\right) dz }{\int_0^1 z(1-z)^{-3/2} e^{\frac{x^2}{2}z} \pi\left( \frac{z}{1-z}\right) dz } 
= \frac{\int_0^\infty u(1+u)^{-3/2} e^{\frac{x^2u}{2(1+u)}} \pi\left( u\right) dz }{\int_0^\infty u(1+u)^{-1/2} e^{\frac{x^2u}{2(1+u)}} \pi\left( u\right) dz }
\label{eq:proof:def:mx}
\end{equation}
Following the proof of lemma A.6. of \cite{vanderpas2016}, we have that working under Conditions \ref{cond1} and \ref{cond2},
\begin{align*}
m_x &= \frac{\int_0^\infty u(1+u)^{-3/2} e^{\frac{x^2u}{2(1+u)}} \pi\left( u\right) dz }{\int_0^\infty u(1+u)^{-1/2} e^{\frac{x^2u}{2(1+u)}} \pi\left( u\right) dz } \\ 
& \leq s_n + \frac{\sqrt{2}}{c} \int_{s_n}^\infty \frac{u}{(1+u)^{3/2}} e^{\frac{x^2}{2} \frac{u}{1+u}} \pi(u) du \\ 
&\leq s_n \left( 1 + \frac{\sqrt{2}}{c} C e^{x^2/4} \right) + \frac{\sqrt{2}}{c} \int_1^{\infty} \frac{u}{(1+u)^{3/2}} e^{\frac{x^2}{2} \frac{u}{1+u}} \pi(u) du \\ 
& \leq s_n \left( 1 + \frac{\sqrt{2}}{c} C e^{x^2/4} \right) + e^{\frac{x^2}{2}}\frac{\sqrt{2}}{c} \int_1^{\infty}  \frac{\pi(u)}{\sqrt{u}} du \\ 
& \leq s_n \left( 1 + \frac{\sqrt{2}}{c} C e^{x^2/4} \right) + e^{\frac{x^2}{2}}\frac{\sqrt{2}}{c} 2Cs_n \frac{1}{\sqrt{\log(n/p_n)}},
\end{align*}
where the last line is obtained by splitting the integral $\int_1^\infty$ in $\int_1^{\log(n/p_n)}$ and $\int_{\log(n/p_n)}^\infty$ and using condition \ref{cond3} on both. %
Using this inequality, we can thus control the type I error of the testing procedure. We have
\begin{align*}
\P_0(m_X \geq \alpha) \leq& \P_0\left(s_n \left( 1 + \frac{\sqrt{2}}{c} C e^{X^2/4} \right) + e^{\frac{X^2}{2}}\frac{\sqrt{2}}{c} 2Cs_n \frac{1}{\sqrt{\log(n/p_n)}} \geq \alpha \right) \\ 
\leq& \P_0\left(s_n \left( 1 + \frac{\sqrt{2}}{c} C e^{X^2/4} \right) \geq \alpha \right) +  \\ 
& \P_0\left(e^{\frac{X^2}{2}}\frac{\sqrt{2}}{c} 2Cs_n \frac{1}{\sqrt{\log(n/p_n)}} \geq \alpha  \right) \\
=& P_1 + P_2
\end{align*}
We will control each term separately. The first term can be bounded by
\begin{align*}
P_1 &= \P_0\left(X^2 \geq 4 \log\left( \frac{c}{Cs_n \sqrt{2}}(\alpha - s_n) \right)  \right) \\ 
&=  4\sqrt{2 \pi} \left(\log\left( \frac{c}{Cs_n \sqrt{2}}(\alpha - s_n) \right) \right)^{-1/2} e^{-2 \log\left( \frac{c}{Cs_n \sqrt{2}}(\alpha - s_n) \right) }(1+o(1)) \\
&\lesssim C \frac{p_n^2}{n^2} \log(n/p_n)^2(1 + o(1)),
\end{align*}
where we have used the fact that $\int_z^\infty x^{s-1} e^{-x/a} dx = az^{s-1}e^{-z/a}(1+o(1))$ as $z\to \infty$. 
We now bound the second term using the same arguments  
\begin{align*}
P_2 &=\P_0 \left(X^2 \geq  2 \log\left( \frac{\alpha c}{2\sqrt{2} C} \frac{\sqrt{\log(n/p_n)}}{s_n}  \right)  \right)  \\ 
& = 2 \sqrt{2 \pi} \log\left( \frac{\alpha c}{2\sqrt{2} C} \frac{\sqrt{\log(n/p_n)}}{s_n}  \right)^{-1/2} e^{-\log\left( \frac{\alpha c}{2\sqrt{2} C} \frac{\sqrt{\log(n/p_n)}}{s_n}  \right)}(1 + o(1)) \\ 
&= 8 \sqrt{\pi} \frac{C}{\alpha c} \frac{p_n}{n} (1+o(1)).
\end{align*}
We thus have that 
\begin{equation}
\P_0(m_X \geq \alpha) \leq 8 \sqrt{\pi} \frac{C}{\alpha c} \frac{p_n}{n}(1 + o(1))
\label{eq:proof:BayesRisk:t1}
\end{equation}

\subsubsection{Control of the Type II Error}

We now control the type II error for both risks i.e. either $\P_{\mathcal{N}(0,\psi)}(\xi = 0)$ or $\P_{\mathcal{N}(\theta,1)}(\xi = 0)$ depending on the considered risk, when the prior satisfies condition \ref{cond1}. 
The following Lemma gives a lower bound on $m_x$ for $x$ large enough, which will prove useful to bound the type two error of each individual tests. 

\begin{lem}
Let $m_x$ be defined as \eqref{eq:proof:def:mx}.
If $\pi$ satisfies condition \ref{cond1}, there exists a constant $C_1$ depending only on $\alpha$, $u_*$ and $u_0$ such that for all $|x| \geq C_1 + \sqrt{2K(1+u_0)\log(n/p)}$ we have
$
m_x \geq \alpha.
$
\label{lem:proof:bound:mxCondition1}
\end{lem}
The proof of this lemma is postponed to appendix \ref{sec:appendix:proof:lem:bound:mxCondition1}. Given this Lemma, we bound the individual type two error, let $t_{2,i} = \P_{1+\psi}(\xi_i = 0 )$ and define $T_n:= C_1 + \sqrt{2K(1+u_0)\log(n/p_n)}$ we have, given the preceding lemma
\begin{align*}
t_{2,i} &= \P_{1+\psi_n}(m_X \leq \alpha \cap X \leq T_n) \\ 
&\leq \P_{1+\psi_n}(X \leq T_n) \\ 
&= 2\Phi\left( \frac{T_n}{\sqrt{1+\psi_n^2}} \right) -1 \\ 
&= \left(2\Phi\left( \sqrt{2K(u_0 + 1)C_\psi} \right) -1\right)(1+o(1)).
\end{align*}
This conclude the proof of Theorem \ref{thm:BayesRisk}. 
\subsection{Proof of Theorem \ref{thm:seprate}}

\subsubsection{False discovery rate}

Recall that for a collection of test $\xi_1, \dots \xi_n$, and a true parameter $\theta$ with support $S_0$ the False Discovery Rate is defined by 
$$
FDR_n(\theta) = \E_\theta \left( \frac{\sum_{i \in S_0^c} \xi_i }{\max(\sum_{i=1}^n \xi_i,1)} \right).
$$
First note that using Lemma \ref{lem:proof:bound:mxCondition1}, under the assumption that $\min_{i\in S_0} |\theta_i| >T_n + v_n$ we have 
$$
t_{2,i} \leq \Phi(T_n - \theta_i) - \Phi(-T_n - \theta_i) \leq  \Phi\left(T_n - \min_{i \in S_0}(|\theta_i|)\right) = o(1)
$$
We now show that the number of true discoveries $\sum_{i \in S_0} \xi_i$ is at least a fraction of $p_n$ with probability that goes to $1$. Let $a_n = \Phi(T_n - \min_{i\in S_0}|\theta_i|)$, and let $t_n = \log(1/a_n)$, then using Chernoff inequality we have

\begin{align*}
\P_\theta(\sum_{i\in S_0} \xi_i &\leq \lambda p_n ) = e^{\lambda t_n p_n}\E(e^{-t_n\sum_{i\in S_0} \xi_i }) \\ 
& \leq e^{-(1-\lambda) t_n p_n + \sum_{i \in S_0} \log(1 + t_{2,i}(e^{t_n} - 1))} \\ 
& \leq e^{-p_n((1-\lambda) t_n  - \log(2 - a_n))} = o(1)
\end{align*}


We thus have the following upper bound for the false discovery rate provided that $\min_{i \in S_0}(|\theta_i|)> T_n + v_n$ 
\begin{align}
FDR_n(\theta) &= \E_\theta \left( \frac{\sum_{i \in S_0^c} \xi_i }{\max(\sum_{i=1}^n \xi_i,1)} \1_{\sum_{i=1}^n \xi_i > \lambda p_n} \right) + o(1) \nonumber \\ 
& \leq \frac{(n-p_n)t_1}{(n-p_n)t_1 + \lambda p_n} + o(1)\label{eq:boundfdr:jensen} \\ 
& \leq \frac{8 \sqrt{\pi} \frac{C}{\alpha c} (1 + o(1))}{(8 \sqrt{\pi} \frac{C}{\alpha c} + \lambda)(1 + o(1))} + o(1)\nonumber \\ 
& = \frac{1}{1 + \frac{\lambda{\alpha c}}{8 \sqrt{\pi} C} (1 + o(1))}+ o(1) \nonumber
\end{align}
where \eqref{eq:boundfdr:jensen} is obtained using Jensen inequality given that the function $x \mapsto \frac{x}{x+a}$ is concave. Thus by choosing $C$ small enough, a bound on the FDR.

\subsubsection{False Nondiscovery rate}

Using the same method, we also bound the False Nondiscovery Rate (FNR) defined as 
$$
FNR_n(\theta) = \frac{\E_\theta\left( \sum_{i\in S_0} (1-\xi_i) \right)}{p_n}
$$
Using the preceding results we have 
\begin{align*}
FNR_n(\theta) &= \frac{1}{p_n}\sum_{i\in S_0} t_{2,i} \\ 
&\leq \Phi(T_n - \min_{i\in S_0} |\theta_i|) + o(1)
\end{align*}
which ends the proof.

\section{Discussion}

We expanded the class of Gaussian scale mixture priors with theoretical guaranties for the multiple testing problems for two different types of risks. The conditions are similar to the ones considered in \cite{vanderpas2016} which indicates that theses priors are well suited both for estimation and testing. What seems to be key in the bounds obtained is the tails of the prior distribution. More precisely, the tails of the local variance should be at least heavy as an exponential to allow for the detection of signals. 

In regards of the results in \cite{Ghosh2015} it seems that the bounds obtained for the Bayes Risk could be sharpen. However, this will require some more restrictive conditions on the prior distribution. To our best knowledge, there were no existing bounds on the minimax risk for any Bayesian multiple testing procedure. The bounds proposed here are competitive with the one found in \cite{arias-castro2017} or \cite{Rabinovich17}. Although not aiming at the same problem \cite{vanderpas2017UQHS} showed that honest uncertainty quantification cannot be attained for true signals bellow a similar threshold, we thus believe that theses bounds are of the correct order.

We only considered here the simplistic Gaussian sequence model, but the results obtained for the Gaussian Scale mixture priors are encouraging and we could thus hope to extend it to more complex high dimensional models. This would be of great interest since these priors are computationally attractive \citep[see][]{MR3620452}, contrariwise to the Spike and Slab priors.

\appendix

\section{Technical Lemmas}

\begin{lem}
\label{lem:upperbound:Lau}
Let $L$ be a uniformly regularly varying function at infinity with constants $u_0$ and $R$. For all $u>u_0$ and $a>1$ we have 
$$
\frac{L(u)}{L(au)} \leq (2a)^{\log_2(R)}
$$
\end{lem}
\begin{lem}
\label{lem:shiftedLtilda}
Let $L$ be a uniformly regularly varying function at infinity with constants $u_0$ and $R$. Then $\tilde{L}(\cdot) = L(\cdot-1)$ is also uniformly regularly varying at infinity with constants $u_0+1$ and $R^3$. 
\end{lem}

\section{Proof of Lemma \ref{lem:proof:bound:mxCondition1}}
\label{sec:appendix:proof:lem:bound:mxCondition1}
\begin{proof}
The proof is adapted from the the proof of Theorem 2.1 in \cite{vanderpas2016}. Given the definition of $m_x$ we have 
$$
1 - m_x = \frac{ \int_1^\infty e^{-\frac{x^2}{2v}} v^{-3/2} \pi(v-1)dv }{ \int_1^\infty e^{-\frac{x^2}{2v}} v^{-1/2} \pi(v-1)dv } = \frac{N_n}{D_n}. 
$$
Denote by $I(a,b) = \int_a^b e^{-\frac{x^2}{2v}} v^{-3/2} \pi(v-1)dv/D_n$, we immediately have 
$$
I\left( \frac{1+u_0}{1-\alpha}, + \infty \right) \leq \frac{1-\alpha}{1+u_0}.
$$
We now control the remaining term $I\left( 1,\frac{1+u_0}{1-\alpha} \right)$. Given that $\pi$ satisfies condition \ref{cond1} and taking $C_1>2u_0$, we have the following lower bound on $D_n$ 
\begin{align}
D_n &\geq \int_{|x|/2}^{|x|} e^{-\frac{x^2}{2v}} v^{-1/2} \pi(v-1)dv \nonumber \\ 
& \geq e^{b - |x|(1+b)} |x|^{1/2} \int_{1/2}^1 L_n(|x|z-1) dz \nonumber \\ 
& \geq \frac{1}{2R^3} e^{b - |x|(1+b)} |x|^{1/2} L_n(|x|/2 - 1), \label{eq:proofs:lowerbound:Dn}
\end{align}
where the last inequality follows from Lemma \ref{lem:shiftedLtilda}. Now consider $I\left( u_0+1,\frac{1+u_0}{1-\alpha} \right)$, we have taking $C_1$ large enough 
\begin{align}
I\left( u_0+1,\frac{1+u_0}{1-\alpha} \right) &= \frac{\int_{u_0+1}^{\frac{u_0+1}{1-\alpha}} e^{-\frac{x^2}{2v}} v^{-3/2} \pi(v-1)dv}{D_n} \nonumber \\ 
&\lesssim e^{-\frac{x^2(1-\alpha)}{2(u_0+1)} + |x|(b+1) } |x|^{-1/2} \int_{u_0+1}^{\frac{u_0+1}{1-\alpha}} \frac{L_n(v-1)}{L_n(\frac{|x|}{2} - 1)} dv \nonumber \\ 
&\lesssim e^{-\frac{x^2(1-\alpha)}{2(u_0+1)} + |x|(b+1)} |x|^{3\log_2(R)-1/2} \label{eq:proof:Imiddle}\\
&\leq \frac{1-\alpha}{2}\frac{u_0}{u_0+1}, \nonumber
\end{align}
where the inequality \eqref{eq:proof:Imiddle} follows from Lemmas \ref{lem:upperbound:Lau} and \ref{lem:shiftedLtilda} and $|x|>C_1$ for $C_1$ a sufficiently large constant that only depends on $R$, $b$, $u_0$ and $\alpha$. We now consider the last term $I(1,u_0+1)$. From equation \eqref{eq:proofs:lowerbound:Dn} we get and condition \ref{cond1} we have 
$$
D_n \geq \frac{1}{2R^3} e^{b + b' - |x|(1+b+b'/2)} |x|^{1/2} \left(\frac{p_n}{n} \right)^K. 
$$
For the numerator we have, given that $\pi$ is a density 
$$
\int_1^{u_0+1} e^{-\frac{x^2}{2v}} v^{-1/2} \pi(v-1)dv \leq e^{-\frac{x^2}{2(u_0+1)}}, 
$$
and thus 
\begin{align*}
I(1,u_0+1) &\leq 2R^3 e^{-\frac{x^2}{2(u_0+1)} + |x|(1+b+b'/2) -b - b' }  |x|^{-1/2} \left( \frac{n}{p_n} \right)^K \\
&\leq \frac{1-\alpha}{2}\frac{u_0}{1+u_0}, 
\end{align*}
given that $|x| \geq u_* + C_1 + \sqrt{2K(1+u_0)\log(n/p_n)}$, where $C_1$ is a constant large enough that only depends on $R$, $b$, $b'$, and $u_0$. This complete the proof .
\end{proof}

\section{Proofs for the adaptive results} 
\label{seq:proof:adaptive}
\subsection{Bayes Risk}

\subsubsection{False Positive}

The number of false positives is given by $\sum_{i=1}^n \1_{\xi_i = 1 \cap v_i = 0} $. 

\begin{align*}
\E\left( \sum_{i=1}^n \1_{\xi_i^{\hat{p}} = 1 \cap v_i = 0}\right) &= \E\left(\1_{\hat{p} \leq C^u p_n} \sum_{i=1}^n \1_{\xi_i^{\hat{p}} = 1 \cap v_i = 0}\right) + o(p_n) \\ 
&\leq  \sum_{i=1}^n \E\left( \1_{\hat{p} \leq C^u p_n} \1_{\xi_i^p = 1} \1_{v_i = 0} \right) + o(p_n)
\end{align*}

Given that $\{ \xi_i^p = 1 \} = \{m_x^p > \alpha\}$ we have that for all $p \leq C^u p_n$ using Conditions 1 and 2 

\begin{align*}
\1\left\{\xi_i^p\right\} &= \1\left\{m_x^p > \alpha\right\} \\ 
& \leq  \1\left\{  \frac{p}{n} \log\left(\frac{p}{n} \right) \left( 1+ \frac{C\sqrt{2}}{c} e^{x^2/4} \right) + \frac{p}{n} \log\left(\frac{p}{n} \right)^{1/2} \frac{2C\sqrt{2}}{c} e^{x^2/2} > \alpha \right\} \\ 
& \leq  \1\left\{  \frac{C^up_n}{n} \log\left(\frac{C^up_n}{n} \right) \left( 1+ \frac{C\sqrt{2}}{c} e^{x^2/4} \right) + \frac{C^up_n}{n} \log\left(\frac{C^up_n}{n} \right)^{1/2} \frac{2C\sqrt{2}}{c} e^{x^2/2} > \alpha \right\}
\end{align*}
We can thus conclude that 
$$
\E\left( \sum_{i=1}^n \1_{\xi_i^{\hat{p}} = 1 \cap v_i = 0}\right) \leq 8 \sqrt{\pi} \frac{C}{\alpha c} C^u p_n(1 + o(1)),
$$
using the same arguments as before. 

\subsubsection{False negative} 

\begin{align*}
\E\left( \sum_{i=1}^n \1_{\xi_i^{\hat{p}} = 0 \cap v_i = 1}\right) &= \E\left(\1_{\hat{p} \geq p_n \left(\frac{p_n}{n}\right)^{-\zeta}e^{-C_d\sqrt{\log(n/p_n)}}} \sum_{i=1}^n \1_{\xi_i^{\hat{p}} = 0 \cap v_i = 1}\right) + o(p_n) \\ 
&\leq  \sum_{i=1}^n \E\left( \1_{\hat{p} \geq c_d p_n \left(\frac{p_n}{n}\right)^{-\zeta}e^{-C_d\sqrt{\log(n/p_n)}}} \1_{\xi_i^{\hat{p}} = 0} \1_{v_i = 1} \right) + o(p_n)
\end{align*}

Let $T_n = C + \sqrt{2K(u_0+1) (1+\zeta)\log(n/p_n)}$ for some $\zeta \geq 0$. Given the preceding results we have if $|x| > T_n$ we have for $C_1$, $C_2$ and $C_3$ some absolute constants 
\begin{align*}
1-m_x^p &\leq \frac{C_1}{|x|}\left( C_2 + e^{(1 + \frac{b + b'}{2})|x| - \frac{x^2}{2u_0} + K\log(n/p) + \frac{1}{2}\log(|x|)} \right) \\ 
& \leq \frac{C_1}{|x|} \left( C_2 +  C_3 \left(\frac{p_n}{p}\right)^K \left(\frac{p_n}{n}\right)^{K\zeta} e^{-KC_d \sqrt{\log(n/p_n)}} \right). 
\end{align*}
For all $|x| > T_n$ 
\begin{align*}
\1\{\xi_i^p = 0\} & = \1\{1-m^p_x \geq  1- \alpha\} \\ 
& \leq \1\left\{ \frac{C_1}{|x|} \left( C_2 +  C_3 \left(\frac{p_n}{p}\right)^K \left(\frac{p_n}{n}\right)^{K\zeta} e^{-KC_d \sqrt{\log(n/p_n)}} \right) > 1-\alpha \right\} 
\end{align*}

Thus 
\begin{align*}
\E\left( \1_{\hat{p} \geq c_d p_n \left(\frac{p_n}{n}\right)^{-\zeta}e^{-C_d\sqrt{\log(n/p_n)}}} \1_{\xi_i^{\hat{p}} = 0} \1_{v_i = 1} \right)  \leq& \E(\1_{|x| \leq T_n} \1_{v_i = 1}) + \\  &\E( \1_{|x|> T_n} \1_{\frac{C_1 C_3 c_d}{|x|}\geq 1-\alpha} \1_{v_i = 1})
\end{align*}
Using the same approach as we did for the deterministic $p$ we get 
$$
\E\left( \sum_{i=1}^n \1_{\xi_i^{\hat{p}} = 0 \cap v_i = 1}\right) \leq p_n\left(2 \Phi\left( \sqrt{2u_0K(1+\zeta)C_\psi} \right) - 1 \right) (1 + o(1)).
$$

\subsection{Minimax risk} 

One can easily adapt the proof of Lemma \ref{lem:proof:bound:mxCondition1} and show that if the prior distribution satisfies condition \ref{cond1} with $p = \hp$ and if $\hp$ satisfies Condition \ref{ass:adapt:Minimax}, then there exists a constant $C_1$ depending only on $u_*$, $u_0$ and $\alpha$ such that for all $|x| \geq C_1 + \sqrt{2K(u_0+1)\log(n/\gamma_n)}$ we have $m_x \leq \alpha$. We prove Theorem \ref{thm:adapt:seprate} following the same lines as for the proof of Theorem \ref{thm:seprate}. 



\bibliographystyle{chicago}
\bibliography{sparseequidae}

\begin{thebibliography}{}

\bibitem[\protect\citeauthoryear{Arias-Castro and Chen}{Arias-Castro and
  Chen}{2017}]{arias-castro2017}
Arias-Castro, E. and S.~Chen (2017).
\newblock Distribution-free multiple testing.
\newblock {\em Electron. J. Statist.\/}~{\em 11\/}(1), 1983--2001.

\bibitem[\protect\citeauthoryear{Baraud}{Baraud}{2002}]{baraud2002}
Baraud, Y. (2002, 10).
\newblock Non-asymptotic minimax rates of testing in signal detection.
\newblock {\em Bernoulli\/}~{\em 8\/}(5), 577--606.

\bibitem[\protect\citeauthoryear{Benjamini and Hochberg}{Benjamini and
  Hochberg}{1995}]{benjamini1995controlling}
Benjamini, Y. and Y.~Hochberg (1995).
\newblock Controlling the false discovery rate: a practical and powerful
  approach to multiple testing.
\newblock {\em Journal of the royal statistical society. Series B
  (Methodological)\/}, 289--300.

\bibitem[\protect\citeauthoryear{Bhattacharya, Chakraborty, and
  Mallick}{Bhattacharya et~al.}{2016}]{MR3620452}
Bhattacharya, A., A.~Chakraborty, and B.~K. Mallick (2016).
\newblock Fast sampling with {G}aussian scale mixture priors in
  high-dimensional regression.
\newblock {\em Biometrika\/}~{\em 103\/}(4), 985--991.

\bibitem[\protect\citeauthoryear{Bogdan, Chakrabarti, Frommlet, and
  Ghosh}{Bogdan et~al.}{2011}]{MR2850212}
Bogdan, M., A.~Chakrabarti, F.~Frommlet, and J.~K. Ghosh (2011).
\newblock Asymptotic {B}ayes-optimality under sparsity of some multiple testing
  procedures.
\newblock {\em Ann. Statist.\/}~{\em 39\/}(3), 1551--1579.

\bibitem[\protect\citeauthoryear{Carvalho, Polson, and Scott}{Carvalho
  et~al.}{2010}]{Carvalho2010}
Carvalho, C.~M., N.~G. Polson, and J.~G. Scott (2010).
\newblock The horseshoe estimator for sparse signals.
\newblock {\em Biometrika\/}~{\em 97\/}(2), 465--480.

\bibitem[\protect\citeauthoryear{Castillo, Schmidt-Hieber, and van~der
  Vaart}{Castillo et~al.}{2015}]{Castillo2015}
Castillo, I., J.~Schmidt-Hieber, and A.~van~der Vaart (2015, 10).
\newblock Bayesian linear regression with sparse priors.
\newblock {\em Ann. Statist.\/}~{\em 43\/}(5), 1986--2018.

\bibitem[\protect\citeauthoryear{Castillo and Van~der Vaart}{Castillo and
  Van~der Vaart}{2012}]{Castillo2012}
Castillo, I. and A.~W. Van~der Vaart (2012).
\newblock Needles and straw in a haystack: Posterior concentration for possibly
  sparse sequences.
\newblock {\em Ann. Statist.\/}~{\em 40\/}(4), 2069--2101.

\bibitem[\protect\citeauthoryear{Efron and Tibshirani}{Efron and
  Tibshirani}{2002}]{efron2002empirical}
Efron, B. and R.~Tibshirani (2002).
\newblock Empirical bayes methods and false discovery rates for microarrays.
\newblock {\em Genetic epidemiology\/}~{\em 23\/}(1), 70--86.

\bibitem[\protect\citeauthoryear{Efron, Tibshirani, Storey, and Tusher}{Efron
  et~al.}{2001}]{MR1946571}
Efron, B., R.~Tibshirani, J.~D. Storey, and V.~Tusher (2001).
\newblock Empirical {B}ayes analysis of a microarray experiment.
\newblock {\em J. Amer. Statist. Assoc.\/}~{\em 96\/}(456), 1151--1160.

\bibitem[\protect\citeauthoryear{Fromont, Lerasle, and Reynaud-Bouret}{Fromont
  et~al.}{2016}]{fromont2016}
Fromont, M., M.~Lerasle, and P.~Reynaud-Bouret (2016, 12).
\newblock Family-wise separation rates for multiple testing.
\newblock {\em Ann. Statist.\/}~{\em 44\/}(6), 2533--2563.

\bibitem[\protect\citeauthoryear{Genovese and Wasserman}{Genovese and
  Wasserman}{2004}]{MR2065197}
Genovese, C. and L.~Wasserman (2004).
\newblock A stochastic process approach to false discovery control.
\newblock {\em Ann. Statist.\/}~{\em 32\/}(3), 1035--1061.

\bibitem[\protect\citeauthoryear{Ghosh and Chakrabarti}{Ghosh and
  Chakrabarti}{2017}]{Ghosh2015}
Ghosh, P. and A.~Chakrabarti (2017).
\newblock Asymptotic optimality of one-group shrinkage priors in sparse
  high-dimensional problems.
\newblock Advance publication.

\bibitem[\protect\citeauthoryear{Ghosh, Tang, Ghosh, and Chakrabarti}{Ghosh
  et~al.}{2016}]{ghosh2016}
Ghosh, P., X.~Tang, M.~Ghosh, and A.~Chakrabarti (2016, 09).
\newblock Asymptotic properties of bayes risk of a general class of shrinkage
  priors in multiple hypothesis testing under sparsity.
\newblock {\em Bayesian Anal.\/}~{\em 11\/}(3), 753--796.

\bibitem[\protect\citeauthoryear{Ingster}{Ingster}{1993}]{Ingster93}
Ingster, Y.~I. (1993).
\newblock Asymptotically minimax hypothesis testing for nonparametric
  alternatives. {I}, {II} and {III}.
\newblock {\em Math. Methods Statist.\/}~{\em 2\/}(2), 85--114.

\bibitem[\protect\citeauthoryear{Johnson and Rossell}{Johnson and
  Rossell}{2012}]{MR2980074}
Johnson, V.~E. and D.~Rossell (2012).
\newblock Bayesian model selection in high-dimensional settings.
\newblock {\em J. Amer. Statist. Assoc.\/}~{\em 107\/}(498), 649--660.

\bibitem[\protect\citeauthoryear{Lehmann}{Lehmann}{1957}]{MR0084952}
Lehmann, E.~L. (1957).
\newblock A theory of some multiple decision problems. {I}.
\newblock {\em Ann. Math. Statist.\/}~{\em 28}, 1--25.

\bibitem[\protect\citeauthoryear{M{\"u}ller, Parmigiani, Robert, and
  Rousseau}{M{\"u}ller et~al.}{2004}]{MR2109489}
M{\"u}ller, P., G.~Parmigiani, C.~Robert, and J.~Rousseau (2004).
\newblock Optimal sample size for multiple testing: the case of gene expression
  microarrays.
\newblock {\em J. Amer. Statist. Assoc.\/}~{\em 99\/}(468), 990--1001.

\bibitem[\protect\citeauthoryear{Polson and Scott}{Polson and
  Scott}{2010}]{Polson2010}
Polson, N.~G. and J.~G. Scott (2010).
\newblock Shrink globally, act locally: Sparse {B}ayesian regularization and
  prediction.
\newblock {\em Bayesian Statistics\/}~{\em 9}, 501--538.

\bibitem[\protect\citeauthoryear{{Rabinovich}, {Ramdas}, {Jordan}, and
  {Wainwright}}{{Rabinovich} et~al.}{2017}]{Rabinovich17}
{Rabinovich}, M., A.~{Ramdas}, M.~I. {Jordan}, and M.~J. {Wainwright} (2017,
  May).
\newblock {Optimal Rates and Tradeoffs in Multiple Testing}.
\newblock {\em ArXiv e-prints\/}.

\bibitem[\protect\citeauthoryear{Robbins}{Robbins}{1956}]{Robbins1956}
Robbins, H. (1956).
\newblock An empirical {B}ayes approach to statistics.
\newblock In {\em Proceedings of the Third Berkeley Symposium on Mathematical
  Statistics and Probability, Volume 1: Contributions to the Theory of
  Statistics}, Berkeley, California, pp.\  157--163. University of California
  Press.

\bibitem[\protect\citeauthoryear{Rossell and Telesca}{Rossell and
  Telesca}{2015}]{rossell2015non}
Rossell, D. and D.~Telesca (2015).
\newblock Non-local priors for high-dimensional estimation.
\newblock {\em Journal of the American Statistical Association\/}.

\bibitem[\protect\citeauthoryear{{Ro\v ckov\'a} and George}{{Ro\v ckov\'a} and
  George}{2016}]{Rockova2015}
{Ro\v ckov\'a}, V. and E.~I. George (2016).
\newblock The spike-and-slab lasso.
\newblock {\em Journal of the American Statistical Association\/}~{\em
  0\/}(ja), 0--0.

\bibitem[\protect\citeauthoryear{Salomond}{Salomond}{2017}]{salomondTest}
Salomond, J.-B. (2017).
\newblock Testing un-separated hypotheses by estimating a distance.
\newblock {\em Bayesian Anal.\/}.
\newblock Advance publication.

\bibitem[\protect\citeauthoryear{Scott and Berger}{Scott and
  Berger}{2006}]{MR2235051}
Scott, J.~G. and J.~O. Berger (2006).
\newblock An exploration of aspects of {B}ayesian multiple testing.
\newblock {\em J. Statist. Plann. Inference\/}~{\em 136\/}(7), 2144--2162.

\bibitem[\protect\citeauthoryear{Tibshirani}{Tibshirani}{1996}]{Tibshirani1996}
Tibshirani, R. (1996).
\newblock Regression shrinkage and selection via the lasso.
\newblock {\em J. R. Stat. Soc. Ser. B Stat. Methodol.\/}~{\em 58\/}(1),
  267--288.

\bibitem[\protect\citeauthoryear{van~der Pas, Kleijn, and van~der
  Vaart}{van~der Pas et~al.}{2014}]{vanderPas2014}
van~der Pas, S., B.~Kleijn, and A.~van~der Vaart (2014).
\newblock The horseshoe estimator: Posterior concentration around nearly black
  vectors.
\newblock {\em Electron. J. Stat.\/}~{\em 8}, 2585--2618.

\bibitem[\protect\citeauthoryear{van~der Pas, Salomond, and
  Schmidt-Hieber}{van~der Pas et~al.}{2016}]{vanderpas2016}
van~der Pas, S., J.-B. Salomond, and J.~Schmidt-Hieber (2016).
\newblock Conditions for posterior contraction in the sparse normal means
  problem.
\newblock {\em Electron. J. Statist.\/}~{\em 10\/}(1), 976--1000.

\bibitem[\protect\citeauthoryear{van~der Pas, Szabó, and van~der
  Vaart}{van~der Pas et~al.}{2017a}]{vanderpas2017}
van~der Pas, S., B.~Szabó, and A.~van~der Vaart (2017a).
\newblock Adaptive posterior contraction rates for the horseshoe.
\newblock {\em Electron. J. Statist.\/}~{\em 11\/}(2), 3196--3225.

\bibitem[\protect\citeauthoryear{van~der Pas, Szabó, and van~der
  Vaart}{van~der Pas et~al.}{2017b}]{vanderpas2017UQHS}
van~der Pas, S., B.~Szabó, and A.~van~der Vaart (2017b).
\newblock Uncertainty quantification for the horseshoe.
\newblock Advance publication.

\end{thebibliography}

\end{document}